\documentclass{elsart}
\usepackage{amssymb, amsmath, amsfonts}
\usepackage{bm}
\usepackage[all,2cell]{xy}\UseAllTwocells

\def\CC{{\mathcal C}}
\def\DD{{\mathcal D}}
\def\EE{{\mathcal E}}

\def\MM{{\mathcal M}}

\def\PP{{\mathcal P}}

\def\tto{\Rightarrow}
\def\pr{\mathrm{pr}}

\def\entw{\mathbf{entw}}
\def\coring{\mathbf{coring}}
\def\ots{\otimes_S}
\def\otr{\otimes_R}
\def\otu{\otimes_U}
\def\otv{\otimes_V}
\def\ota{\otimes_A}
\def\otb{\otimes_B}
\def\ote{\otimes_E}
\def\comc{\mathrm{\bf Comc}} 
\def\bimod{\operatorname{\mathrm{Bimod}}}
\def\catbimod{\mathbf{bimod}}

\def\nxpoint{\refstepcounter{subsection}%
  \makepoint{\thesubsection}}

\def\makepoint#1{\medbreak\noindent{\bf #1. }}
\newcommand{\nodo}[1]{}

\begin{document}\begin{frontmatter}

\title{Bicategory of entwinings}

\author{Zoran \v{S}koda}
\address{Theoretical Physics Division,
Institute Rudjer Bo\v{s}kovi\'{c}, Bijeni\v{c}ka cesta~54, P.O.Box 180,
HR-10002 Zagreb, Croatia}
\ead{zskoda@irb.hr}

\date{{}}

\maketitle
\begin{abstract}
\noindent 
     We define a bicategory in which the 0-cells are the
entwinings over variable rings. The 1-cells are triples of a
bimodule and two maps of bimodules which satisfy an additional
hexagon, two pentagons and two (co)unit triangles; and the 2-cells
are the maps of bimodules satisfying two simple compatibilities.
The operation of getting the ``composed coring'' from a given
entwining, is promoted here to a canonical morphism of
bicategories from a bicategory of entwinings to the {\sc Street}'s
bicategory of corings.
\end{abstract}
\end{frontmatter}
\section{Algebras, coalgebras, corings}

\nxpoint The main results of this paper are an explicit
construction of a bicategory whose objects are entwinings and of a
homomorphism of bicategories from that bicategory to {\sc R.
Street}'s bicategory of corings (a straightfoward analogue of the
2-category of (co)monads from \cite{Street:formal}).

\nxpoint The results in this manuscript hold in one of the
following two generalities. In the first case we consider
entwinings $(R,A,C,\psi:C\otimes_R A\to A\otimes_R C)$ between
$R$-algebras $A$ and $R$-coalgebras $C$ over (variable)
commutative unital ring $R$
(\cite{BrzMajid:coalgBundles,brzWisbBk}); in the latter case the
entwinings are between $R$-rings $A$ (monoids in $R-\bimod$) and
$R$-corings $C$ (comonoids in $R-\bimod$) over (variable) not
necessarily commutative ring $R$. I will talk about ``algebras''
and ``coalgebras'' over commutative rings but will be careful
about the sides for the bimodules even when over the ground ring
$R$, $S$ etc., so that all our calculations seem identical in both
cases. $\psi$ does not need to be invertible for our purposes.

One can go to the third generality, working with the internal
entwining structures in monoidal categories
(cf.\cite{Bohm:internal}). In fact, in familiar cases of our
interest, the monoidal categories have obvious coherences for
associativity so we do not write them in our statements and
proofs. On the other hand if the coherences are indeed
algebraically nontrivial, then the statements here are more
complicated and somewhat more interesting.

\nxpoint When the preprint version 1 of this article has been
posted, {\sc G. B\"ohm} has kindly called my attention to the
following argument: distributive laws
(\cite{Semtriples,Beck:distr}) are just monads in a 2-category of
monads in the sense of formal monad theory (\cite{Street:formal})
of {\sc R.~Street}, and in particular they themselves make a
2-category; the analogue can be easily written out for mixed
distributive laws betwen a monad and a comonad; it is not
published in detail, but it is widely known among the experts that
the formal monad theory can be extended to bicategorical setup,
instead of strict 2-categories; finally entwinings are {\it mixed}
distributive laws in the setup of the bicategory of rings and
bimodules; regarding that in a bicategory we can do 2 dualizations
(inverting 1-cells and 2-cells) there are thus 4 natural
bicategories of entwinings. Our construction is explicit and from
scratch and does not use this chain of constructions and
translations of data (explicit descriptions are its merit but also
its conceptual deficiency).

\nxpoint The bicategory introduced here is also analogous to the
2-category of distributive laws between actions of a fixed
monoidal category $\CC$ on variable category $\DD$ and monads in
$\DD$ (such distributive laws were studied in our earlier preprint
\cite{skoda:distr} and the 2-category they form
in~\cite{skoda:ent}); such distributive laws can be identified
with $\CC$-equivariant monads, that is, monads in the 2-category
$\mathbf{act}_c(\CC)$ of $\CC$-actegories, colax $\CC$-equivariant
functors and their natural transformations. Finally, one can study
the distributive laws between actions of two different variable
categories on their common target category. If we restrict in the
latter case to the invertible distributive laws, then we call such
data biactegories (not a typo!); in work in progress
\cite{skoda:biact} we introduce and study a tensor product of
biactegories using certain pseudcoequalizer (the induction
pseudofunctor for actegories can be viewed as a special case of
that pseudocoequalizer) with motivation in associating 2-vector
bundles to principal bundles with structure 2-group. Finally,
biactegories make a tricategory $\mathbf{biact}$ which is a
categorification of the bicategory of bimodules $\catbimod$.

\section{Bicategory $\entw$}

\nxpoint
{\bf 0-cells} of $\entw$ are the entwinings $(R,A,C,\psi)$.
Here ring $R$ may vary!
$\psi:C\otr A\to A\otr C$ satisfies the usual two pentagons
and two triangles (\cite{brzWisbBk}) which we do not write here; these
data are implicit: multiplication $\mu^A: A\otr A\to A$,
unit $\eta^A : R\to A$, comultiplication $\Delta^C : C\to C\otr C$,
counit $\epsilon^C : C\to R$.

\nxpoint
{\bf 1-cells} of $\entw$ are the triples
$(M,\alpha,\beta):(R,A,C,\psi)\to(S,B,D,\chi)$,
where ${}_S M_{R}$ is a $S-R$-bimodule,
and $\alpha : B\ots M\to M\otr A$,
and $\gamma : D\ots M\to M\otr C$ are maps of $S-R$-bimodules,
for which the following 5 diagrams commute:

the hexagon:
\begin{equation}\label{eq:hexagon}
\xymatrix{D\ots B\ots M\ar[d]^{\chi\otimes M}\ar[r]^{D\otimes\alpha} &
D\ots M\otr A\ar[r]^{\gamma\otimes A}&
M\otr C\otr A\ar[d]^{M\otimes\psi}
\\
B\ots D\ots M \ar[r]^{B\otimes\gamma}&
B\ots M\otr C \ar[r]^{\alpha\otimes C} &
M\otr A\otr C
}\end{equation}
the pentagon for the $S-R$-bimodule map $\alpha$:
\begin{equation}\label{eq:alphapentagon}\xymatrix{
B\ots B\ots M\ar[r]^{B\otimes\alpha}\ar[d]^{\mu^N\otimes M}
& B\ots M\otr A\ar[r]^{\alpha\otimes A}
& M\otr A\otr A\ar[d]^{M\otimes\mu^A}\\
B\ots M\ar[rr]^\alpha&& M\otr A
}\end{equation}
the pentagon for the $S-R$-bimodule map $\gamma$:
\begin{equation}\label{eq:gammapentagon}\xymatrix{
D\ots M\ar[rr]^\gamma\ar[d]_{\Delta^D\otimes M}&&
M\otr C\ar[d]^{M\otimes\Delta^C}\\
D\ots D\ots M\ar[r]^{D\otimes\gamma}& D\ots M\otr C\ar[r]^{\gamma\otimes C}&
M\otr C\otr C
}\end{equation}
and the two triangles:
$$\xymatrix{
&M\ar[dl]_{\eta^B\otimes M}\ar[dr]^{M\otimes\eta^A}
&&D\ots M\ar[rr]^\gamma\ar[dr]_{\epsilon^D} &&M\otr C\ar[dl]^{\epsilon^C}\\
B\ots M\ar[rr]^\alpha&& M\otr A&&M
}$$

Notice that in these diagrams we did not bother inserting the brackets and
associativity isomorphisms inherited from the bicategory of bimodules.
A pedantic reader will easily 'correct' this.

\nxpoint {\bf 2-cells} $\theta:(M,\alpha,\gamma)\tto(N,\beta,\delta)$ are the
$S-R$-bimodule maps $\theta:M\to N$ such that
the following two squares commute:
\begin{equation}\label{eq:twosquares}
\xymatrix{
B\ots M\ar[r]^\alpha\ar[d]^{B\otimes\theta}&M\otr A\ar[d]^{\theta\otimes A}
&&
D\ots M\ar[r]^\gamma\ar[d]^{D\otimes\theta}&M\otr C\ar[d]^{\theta\otimes C}
\\
B\ots N\ar[r]^\beta& N\otr A&&D\ots N\ar[r]^\delta& N\otr C
&&
}\end{equation}

\nxpoint {\bf (Composition of 1-cells)}
Given a diagram of morphisms of entwinings
\begin{equation}\label{eq:diagMorEntw}\xymatrix{
(R,A,C,\psi)\ar[r]^{(M,\alpha,\gamma)}&
(S,B,D,\chi)\ar[r]^{(P,\sigma,\tau)}&
(U,E,G,\lambda)\ar[r]^{(Q,\rho,\nu)}&
(V,F,H,\xi)
}\end{equation}
define the composition (up to coherences again)
\begin{equation}\label{eq:compMorEntw}
(P,\sigma,\tau)\circ (M,\alpha,\beta) :=
({}_U P\ots M_R,(P\ots \alpha)\circ(\sigma\ots M),
(P\ots \gamma)\circ(\tau\ots M))
\end{equation}
If we should insert the coherences from the underlying
bicategory of bimodules, instead of $(P\ots\alpha)\circ(\sigma\ots M)$
write the composition of 5 maps
$$\begin{array}{l}
E\otu(P\ots M)\stackrel{a_{E,P,M}}\longrightarrow
(E\otu M)\ots M\stackrel{\sigma\otimes M}\longrightarrow
(P\ots B)\ots M\stackrel{a^{-1}_{P,B,M}}\longrightarrow\\
\stackrel{a^{-1}_{P,B,M}}\longrightarrow
P\ots(B\ots M)\stackrel{P\otimes\alpha}\longrightarrow
P\ots(M\otr A)\stackrel{a_{P,M,A}}\longrightarrow
(P\ots M)\otr A
\end{array}$$
where the isomorphisms
$a_{X,Y,Z} : X\otimes (Y\otimes Z)\to (X\otimes Y)\otimes Z$
are the components of the associtivity coherence of $\catbimod$.

\nxpoint {\bf (Coherences for the composition of 1-cells)}
The coherences from the underlying bicategory $\catbimod$
play the role of coherence in $\entw$ as well, which will be denoted
by the same letter $a$:
\begin{equation}\label{eq:assocCoh}
a_{Q,P,M}:(Q,\mu,\nu)\circ((P,\sigma,\tau)\circ(M,\alpha,\beta))
\to ((Q,\mu,\nu)\circ(P,\sigma,\tau))\circ(M,\alpha,\beta)
\end{equation}
In other words, $a_{(Q,\mu,\nu),(P,\sigma,\tau),(M,\alpha,\beta)}
:= a_{Q,P,M}$ (in particular, the pentagon for the coherence
follows from the pentagon for $a$ in $\catbimod$). We need to
check that this definition is meaningful. The components of $a$
should be 2-cells in $\catbimod$, which are in fact 2-cells in
$\entw$, i.e. two naturality squares commute, first of which is
\begin{equation}\label{eq:triplesq}\xymatrix{
F\otv(Q\otu(P\ots M))\ar[rr]^{``\rho(\sigma\alpha)''}
\ar[d]_{F\otu a_{Q,P,M}}&&
(Q\otu(P\ots M))\otr A\ar[d]^{a_{Q,P,M}\otr A}
\\
F\otv((Q\otu P)\ots M))\ar[rr]^{``(\rho\sigma)\alpha''}&& ((Q\otu
P)\ots M))\otr A, }\end{equation} where the horizontal arrows just
symbollically denoted $''\rho(\sigma\alpha)'',
``(\rho\sigma)\alpha''$ correspond to the analogue of $\alpha$ for
the triple composition, which do involve $\rho,\sigma,\alpha$ but
also various coherences. We will show that this square commutes,
while another square similar to~(\ref{eq:triplesq}), but for the
analogue of $\gamma$ is commutative as well, with similar proof
left to the reader. Keeping track of coherences, the analogue of
$\alpha$ for a single composition is a composition of 5 maps, thus
if for one of them we insert another composition of 5 maps, we
have the composition of 9 maps. In the diagram showing that
(\ref{eq:triplesq}) commutes, and where maps
$''\rho(\sigma\alpha)'', ``(\rho\sigma)\alpha''$ are explicit
compositions of 9 maps each, we will skip tensor product signs,
and enclose the corners of~(\ref{eq:triplesq}) in boxes for
emphasis:
$$\xymatrix{
\fbox{$F(Q(PM))$}\ar[r]^{\,\,\,\,\,\,a_{F,Q,PM}}\ar[d]^{Fa_{Q,P,M}}
&(FQ)(PM)\ar[r]^{\rho (PM)}\ar[d]^{a_{FQ,P,M}}
&(QE)(PM)\ar[r]^{a_{Q,E,PM}}\ar[d]^{a_{QE,P,M}}
&Q(E(PM))\ar[r]^{Q a_{E,P,M}}
&Q((EP)M)\ar[dd]_{Q(\sigma M)}\ar[dl]_{a_{Q,EP,M}}
\\
\fbox{$F((QP)M)$}\ar[d]^{a_{F,QP,M}}
&((FQ)P)M\ar[r]^{(\rho P)M}
&((QE)P)M\ar[r]^{a_{Q,E,P}M}
&(Q(EP))M\ar[d]_{(Q\sigma)M}
&
\\
(F(QP))M\ar[ru]_{a_{F,Q,P}M}&
&&(Q(PB))M\ar[d]_{a_{Q,P,BM}}
&Q((PB)M)\ar[dd]_{Q a_{P,B,M}}\ar[l]_{a_{Q,PB,M}}
\\
\fbox{$((QP)M)A$}
&(QP)(MA)\ar[l]_{\,\,\,\,a_{QP,M,A}}
&(QP)(BM)\ar[l]_{(QP)\alpha}
&((QP)B)M\ar[l]_{a_{QP,B,M}}
&
\\
&\fbox{$(Q(PM))A$}\ar[lu]^{a_{Q,P,M}A}
&Q((PM)A)\ar[l]^{\,\,\,\,a_{Q,PM,A}}
&Q(P(MA))\ar[l]^{Q a_{P,M,A}}\ar[ull]^{a_{Q,P,MA}}
&Q(P(BM))\ar[llu]_{a_{Q,P,BM}}\ar[l]^{Q(P\sigma)}
}$$
The diagram is commutative as it splits into 3 naturality squares
and 4 associativity pentagons. A similar diagram may be written for $\gamma$
instead of $\alpha$.

\nxpoint The composition formula~(\ref{eq:compMorEntw}) indeed
defines a morphism of entwinings. With skipping tensor product signs in
a big diagram I just draw the commutative diagram for showing the hexagon
for the composition.
$$\xymatrix{
GEPM\ar[rr]^{G(P\alpha\circ\sigma M)}\ar[rd]_{G\sigma M}
\ar[ddd]_{\lambda PM}
&& GPMA\ar[rr]^{(P\gamma \circ \tau M)A}\ar[rd]^{\tau MA}
&& PMCA\ar[ddd]^{PM\psi}
\\
& GPBM\ar[ur]^{GP\alpha}\ar[r]_{\tau BM}
& PDBM\ar[r]_{PD\alpha}\ar[d]^{P\chi M}
& PDMA\ar[ur]_{P\gamma A} &\\
& EPDM\ar[dr]_{EP\gamma}\ar[r]^{\sigma DM}
& PBDM\ar[r]^{PD\gamma}
& PBMC\ar[dr]^{P\alpha C}&\\
EGPM\ar[ur]^{E\tau M}\ar[rr]_{E(P\gamma\circ \tau M)}
&&EPMC\ar[ur]_{\sigma MC}\ar[rr]_{(P\alpha\circ\sigma M)C}
&&PMAC
}$$
Horizontal and vertical composition of 2-cells is simply given
by the horizontal and vertical composition of the underlying
morphisms of bimodules. One checks that the composed 2-cells
are indeed 2-cells, i.e. satisfy the two squares~(\ref{eq:twosquares}).

\section{Homomorphism into bicategory of corings}

\nxpoint In this article I work with the old Street's version
$\coring$ of the bicategory of corings. It is summarized in the
end of the article~\cite{brz:bicat} but it is just a variant of a
construction in \cite{Street:formal}). There is a different
(variant of the) bicategory of corings also defined in
\cite{brz:bicat} (could be obtained using certain dualization if
compared to $\coring$) and studied in more detail.\vskip .06in

Let us now define the morphism of bicategories 'composed coring'
$$\fbox{$\comc: \entw\to\coring$}$$

\nxpoint {\bf $\comc$ on objects}:
standard ``composed comonad'' formula well known in coring setup:
$(R,A,C,\psi)$ gives rise to the composed $A$-coring
$(A\otr C,\Delta^\psi,\epsilon^\psi)$ where
$\Delta^\psi$ is the composition
$A\otr C \stackrel{A\otr\Delta^C}\longrightarrow A\otr C\otr C\cong
(A\otr C)\otimes_A (A\otr C)$ and similarly for $\epsilon^\psi$:
$A\otr C\stackrel{A\otr\epsilon^C}\longrightarrow A\otr R\cong A$.

\nxpoint Let us define $\MM\in B-A-Bimod$. the underlying module is
$\MM = M\otr A$, it is a left $B$-module via action
$$
B\ots M\otr A\stackrel{\alpha\otimes A}\longrightarrow
M\otr A\otr A\stackrel{M\otimes\mu^A}\longrightarrow
M\otr A
$$
and a right $A$-module via action $M\otimes\mu^A : M\otr A\otr A\to M\otr A$.
It is easy to check that the two actions are compatible, making $M\otr A$
a $B-A$-bimodule:
$$\xymatrix{
B\ots M\otr A\otr A \ar[d]_{B\otimes M\otimes \mu^A}
\ar[rr]^{\alpha\otimes A\otimes A}
&& M\otimes A\otimes A\otimes A\ar[d]_{M\otimes A\otimes\mu^A}
\ar[rr]^{M\otimes\mu^A\otimes A}
&& M\otr A\otr A\ar[d]_{M\otimes\mu^A}
\\
B\ots M\otr A\ar[rr]^{\alpha\otimes A}
&& M\otr A\otr A\ar[rr]^{M\otimes\mu^A}
&& M\otr A
}$$
Notice that the definition of the $B-A$-module structure on $M\otr A$
implies that the diagram
\begin{equation}\label{eq:maproj}\xymatrix{
B\ots M\otr A\ar[rr]^{\alpha\otimes A}\ar[d]^{\pr}&&
M\otr A\otr A\ar[d]^{\pr}
\\
B\otb(M\otr A)\ar[r]^{\cong}& M\otr A\ar[r]^\cong&(M\otr A)\ota A
}\end{equation}
commutes.

\nxpoint {\bf $\comc$ on morphisms}: the triple
$({}_S M_R,\alpha,\gamma) : (R,A,C,\psi)\to(S,B,D,\chi)$
maps to the pair $(\MM,\zeta)$ where $\MM = M\otr A$
is a $B-A$-bimodule as above and for the map $\zeta$,
one first defines an auxiliary map
$\bar\zeta :B\ots D\ots (M\otr A)\to (M\otr A)\otr A\otr C$
as the composition
$$
B\ots D\ots M\otr A\stackrel{B\otimes\gamma\otimes A}\longrightarrow
B\ots M\otr C\otr A\stackrel{\alpha\otimes\psi}\longrightarrow
M\otr A\otr A\otr C,
$$
One checks that $\bar\zeta$ is a map of $B-A$-bimodules.
The fact that $\bar\zeta$ is a map of left $B$-modules essentially boils
down to the pentagon for map $\phi$. I will skip the tensor signs
in drawing the commutative diagram showing this:
$$\xymatrix{
BBDMA\ar[r]^{BB\gamma A} \ar[dd]^{\mu^B DMA}
& BBMCA\ar[r]^{B\gamma CA}\ar[dd]^{\mu^B MCA}
& BMACA\ar[r]^{BMA\psi}\ar[d]_{\alpha ACA}
& BMAAC\ar[d]^{\alpha AAC}
\\
&& MAACA\ar[d]_{M\mu^A CA}\ar[r]^{MAA\psi}
& MAAAC\ar[d]^{M\mu^A AC}
\\
BDMA\ar[r]_{B\gamma A}
& BMCA \ar[r]_{\alpha CA}
& MACA\ar[r]_{MA\psi}
& MAAC
}$$
Similarly, the fact that $\bar\zeta$ is a map of right $A$-modules
similarly essentially boils to the pentagon for entwining $\psi$:
$$\xymatrix{
BDMAA\ar[r]^{B\gamma AA}\ar[dd]^{BDM\mu^A}
& BMCAA\ar[r]^{\alpha CAA}\ar[dd]^{BMC\mu^A}
& MACAA\ar[r]^{MA\psi A}\ar[dd]^{MAC\mu^A}
& MAACA\ar[d]^{MAA\psi}
\\
&&&MAAAC\ar[d]^{MA\mu^A C}\\
BDMA\ar[r]^{B\gamma A}
&BMCA\ar[r]^{\alpha CA}
&MACA\ar[r]^{MA\psi}
& MAAC
}$$
Let $\nu_1 :(B\ots D)\ots(M\otr A)\to (B\ots D)\otb(M\otr A)$
and $\nu_2 : M\otr A\otr A\otr C\to M\otr A\ota A\otr C \cong M\otr A\otr C$
be the canonical projections.
Now once we defined $\bar\zeta : B\ots D\ots M\otr A\to M\otr A\otr A\otr C$
we prove that there is a unique map
$\bar\zeta':(B\ots D)\otb(M\otr A)\to M\otr A\otr A\otr C$ such that
if $\nu_1\circ\bar\zeta' = \bar\zeta$ and then define
$\zeta:=\nu_2\circ\bar\zeta'$. Showing this is a longer naturality
calculation, involving the hexagon for the map $(M,\alpha,\gamma)$,
the pentagon for $\alpha$, the pentagon for the entwining $\psi$
and 4 or 6 naturality squares (depending on the way of defining $\bar\zeta$).
Indeed, start with $(B\ots D)\ots B\ots (M\otr A)$ and acts with middle $B$
either to the first or the second tensored pair. After that
apply $\bar\zeta$.
Thus, omiting the tensor (over the ground rings)
signs, the composition
$$
BDBMA\stackrel{BD\alpha A}\longrightarrow BDMAA
\stackrel{BDM\mu^A}\longrightarrow
BDMA\stackrel{B\gamma A}\longrightarrow
BMCA\stackrel{BM\psi}\longrightarrow
BMAC\stackrel{\alpha AC}\longrightarrow
MAAC
$$
equals the composition
$$
BDBMA\stackrel{B\chi MA}\longrightarrow
BBDMA\stackrel{\mu^B DMA}\longrightarrow
BDMA\stackrel{B\gamma A}\longrightarrow
BMCA\stackrel{BM\psi}\longrightarrow
BMAC\stackrel{\alpha AC}\longrightarrow
MAAC
$$
(Warning: if one parallely truncates the tail of the two chains of maps, one
does not get an identity!)
$$\xymatrix{
BDBMA\ar[r]^{BD\alpha A}\ar[d]^{B\chi MA}
& BDMAA\ar[r]^{B\gamma AA}
& BMCAA\ar[r]^{BMC\mu^A}\ar[dr]^{\alpha CAA}\ar[d]_{BM\psi A}
& BMCA\ar[r]^{BM\psi}\ar[dr]^{\alpha CA}
& BMAC\ar@/^2pc/[ddd]_{\alpha AC}
\\
BBDMA\ar[r]^{BB\gamma A}\ar[dd]^{\mu^B DMA}
& BBMCA\ar[r]^{B\alpha CA}\ar[dd]^{\mu^B MCA}
& BMACA\ar[d]_{\alpha ACA}
& MACAA\ar[r]_{MAC\mu^A}\ar[ld]^{MA\psi A}
& MACA\ar[dd]_{MA\psi}
\\ &&MAACA\ar[d]_{M\mu^A CA}\ar[r]_{MAA\psi}&MAAAC\ar[dr]_{M\mu^A AC}&
\\ BDMA\ar[r]_{B\gamma A}
&BMCA\ar[r]_{\alpha CA}
&MACA\ar[rr]_{MA\psi}
&&MAAC
}$$

\nxpoint In the situation
\begin{equation}\label{eq:diagMorEntw2}\xymatrix{
(R,A,C,\psi)\ar[rr]^{(M,\alpha,\gamma)}&&
(S,B,D,\chi)\ar[rr]^{(P,\sigma,\tau)}&&
(U,E,G,\lambda)
}\end{equation}
consider the diagram of corings
\begin{equation}\label{eq:diagMorCor}\xymatrix{
\CC\ar[rr]^{(\MM,\zeta^M)}&&
\DD\ar[rr]^{(\PP,\zeta^P)}&&
\EE
}\end{equation}
where the corings are
$\CC = (A\otr C,\Delta^\psi,\epsilon^\psi)$,
$\DD = (B\ots D,\Delta^\chi,\epsilon^\chi)$,
$\EE = (E\otu G,\Delta^\lambda,\epsilon^\lambda)$,
(over $A,B$ and $E$, respectively)
and $\MM = {}_B\MM_A= M\otr A$ and $\PP = {}_E\PP_B= P\ots B$
are the corresponding bimodules.

\nxpoint {\bf Proposition.} {\it The pair $({}_B \MM_{A}, \zeta)
= (M\otr A,\zeta)$ defined above
is a 1-cell in Street's $\coring$. In other words, the pentagon
\begin{equation}\label{eq:streetonecell}\xymatrix{
\DD\otb\MM\ar[rr]^\zeta\ar[d]_{\Delta^\CC\otb\MM}
&&\MM\ota\CC\ar[d]^{\MM\ota\CC}
\\
\DD\otb\DD\otb\MM\ar[r]^{\DD\otb\zeta}
&\DD\otb\MM\ota\CC\ar[r]^{\zeta\ota\Delta^\CC}
&\MM\ota\MM\ota\CC
}\end{equation}
commutes and the compatibility with the counits holds.
}

{\it Proof.} In fact, we shall prove the commutativity of a
diagram in which the upper row is the representative of the map
$\zeta$ at the level of the tensor products over $S$ and $R$, but
in the lower row we indeed have the equivalence classes. The
diagrams are a bit more complicated so we will in addition to
skipping the tensor products over $S$ and $R$, also abbreviate the
signs for the tensor products over $A$ and $B$ by a dot $\cdot$
(the modules involved gurantee that the meaningful choice between
$\ota$ and $\otb$ is unique). For example, $BD\cdot MA$ means
$(B\ots D)\otb(M\otr A)$. One also needs to be careful
``cancelling'' $B$ and $A$ in tensor products over $B$ and $A$
respectively. Carefully distinguish the following two maps (and
their analogues). The first is the natural projection from the
tensor product over $S$ to the tensor product over $M$, say
$BDBM\stackrel{BD\cdot BM}\longrightarrow BD\cdot BM$ (also
sometimes shortly denoted $\mathrm{pr}$) and another is inserting
the unit over $B$, say the map $BD\cdot \eta^B M : BDM\cong
BDSM\to BD\cdot BM$. Thus diagram ~(\ref{eq:streetonecell}) may be
expanded to
$$\xymatrix{
BDMA\ar[rr]^{B\gamma A}\ar[d]^{B\Delta^D MA}
&&BMCA\ar[r]^{BM\psi}\ar[d]_{BM\Delta^C A}
&BMAC\ar[r]^{\alpha AC}\ar[rd]^{BMA\Delta^C}
&MAAC\ar@/^2pc/[dd]^{MA\Delta^C}
\\
BDDMA\ar[r]^{BD\gamma A}\ar[d]^{BD\cdot \eta^B DMA}
&BDMCA\ar[r]^{B\gamma CA}\ar[d]^{BD\cdot \eta^B MCA}\ar[dr]^{BDM\psi}
&BMCCA\ar[r]^{BMC\psi}
&BMCAC\ar[r]^{BMC\psi}
&BMACC\ar[d]_{\alpha AC}
&
\\
BD\cdot BDMA\ar[r]^{BD\cdot B\gamma A}\ar[dd]^\pr
& BD\cdot BMCA\ar[d]_{BD\cdot BM\psi}
& BDMAC\ar[d]^{BD\cdot M\eta^A AC}\ar[dl]^{BD\cdot\eta^B MAC}
\ar[ru]^{B\gamma AC}\ar@/^2pc/[rd]^{\pr}
&&MAACC\ar[d]_{\pr}&
\\
&BD\cdot BMAC\ar[r]_{BD\cdot\alpha AC}
&BD\cdot MAAC\ar[d]_\pr
&BD\cdot MAC\ar[r]^{\zeta C}\ar[dl]^\pr
& MA\cdot ACC\ar[d]_\pr &
\\
BD\cdot BD\cdot MA\ar[rr]^{BD\cdot\zeta} && BD\cdot MA\cdot
AC\ar[rr]^{\zeta\cdot AC} && MA\cdot AC\cdot AC & }$$ where one
can directly observe the commutativity of all smallest circuits
and hence of the entire diagram.

In addition to the pentagon~(\ref{eq:streetonecell}), we need to
check the compatibility of the map $\zeta$ with the counits of the
corings involved: $\epsilon^\DD\otb (M\otr A) = ((B\ots
M)\ota\epsilon^\CC)\circ\zeta$. This follows by the calculation
for the representatives, namely the diagram
$$\xymatrix{
B\ots D\ots M\otr A\ar[r]^{B\otimes\gamma\otimes A}
\ar[d]_{B\otimes\epsilon^D\otimes M\otimes A}&
B\ots M\otr C\otr A\ar[r]^{\alpha\otimes C\otimes A}
\ar[dl]_{B\otimes M\otimes\epsilon^C\otimes A}&
M\otr A\otr C\otr A\ar[d]^{M\otimes A\otimes \psi}
\ar[dl]_{M\otimes A\otimes\epsilon^C\otimes A}&
\\
B\ots M\otr A\ar[r]^{\alpha\otimes M\otimes A}\ar[d]^\pr
& M\otr A\otr A\ar[dr]^\pr &M\otr A\otr A\otr C
\ar[l]_{M\otimes A\otimes A\otimes\epsilon^C} \\
B\otb (M\otr A)\ar[r]^{\cong}& M\otr A&(M\otr A)\ota A\ar[l]^\cong&
}$$
commutes. This finishes the proof.

\nxpoint {\bf $\comc$ on 2-cells}:
If $\theta : (M,\alpha,\gamma)\to (N,\beta,\delta)$ then
$(\theta:M\to N)\mapsto(\theta\otr A:M\otr A\to N\otr A)$. One sees
easily that $\theta\otr A$ is indeed a map of $B-A$-bimodules.
We will just draw the diagram for the left $B$-equivariance:
$$\xymatrix{
B\ots M\otr A\ar[d]_{B\otimes\theta\otimes A}\ar[rr]^{\alpha\otimes A}
&& M\otr A\otr A\ar[d]_{\theta\otimes A\otimes A}\ar[rr]^{M\otimes\mu^A}
&&M\otr A\ar[d]^{\theta\otimes A}\\
B\ots N\otr A\ar[rr]^{\alpha\otimes A}
&& N\otr A \otr A\ar[rr]^{N\otimes\mu^A}
&& N\otr A
}$$
For a fixed domain and codomain 1-cells,
this tautological map is injective, but not the surjective map,
because the property that the composition is a 2-cell in $\coring$ is
weaker than the property that $\theta\otr A$ is a 2-cell in $\entw$. In
$\coring$ case just the external square in
$$\xymatrix{
B\ots D\ots M\otr A\ar[r]^{B\otimes\gamma\otimes A}
\ar[d]_{B\otimes D\otimes\theta\otimes A} &
B\ots M\otr C\otr A\ar[r]^{\alpha\otimes\psi}
\ar[d]_{B\otimes \theta\otimes\psi}&
M\otr A\otr A\otr C\ar[d]^{\theta\otimes A\otimes A\otimes C}\\
B\ots D\ots N\otr A\ar[r]_{B\otimes\delta\otimes A}&
B\ots N\otr C\otr A\ar[r]_{\gamma\otimes\psi}&
N\otr A\otr A\otr C
}$$
commutes. The commutativity of the
right-hand square is implied form one of the squares in the axioms
for $\theta$ ($\theta$ vs. $\gamma$),
while the left-hand square is actually the pasting of
another such square ($\theta$ vs. $\alpha$) and of a naturality
square for the tensoring with $\psi$.

\nxpoint We need to check the functoriality.
Thus consider again the diagram
~(\ref{eq:diagMorEntw}) of morphisms of entwinings and the
compositions~(\ref{eq:compMorEntw}).
$$
\comc(P,\sigma,\tau)\circ_\coring\comc(M,\alpha,\beta)
= (\PP\otb\MM,(\zeta^P\otb\MM)\circ(\PP\ota\zeta^M)),
$$
versus
$$
\comc(P\ots M, (P\ots\alpha)\circ(\sigma\ots M),(P\ots\gamma)\circ
(\tau\ots M))=((P\ots M)\otr A, \zeta^{P\otimes M})
$$
Up to coherences (some of which we already skipped), the two answers
should agree; the additional coherences make the functoriality
true up to invertible 2-cell (pseudofunctoriality).
If we look at the underlying module, this is obvious:
$(P\ots M)\otr A\cong (P\ots B)\otb (M\otr A) = \PP\otb\MM$,
and the agreement for $\zeta$-s is
$$\xymatrix{
\EE\ote(\PP\otb\MM)\ar[rrr]^{\zeta^{P\otimes M}}\ar[d]^\cong &&&
(\PP\otb\MM)\ota\CC\ar[d]^\cong\\
(\EE\ote\PP)\otb\MM\ar[r]^{\zeta^P\otb\MM}
&(\PP\otb\DD)\otb\MM\ar[r]^\cong
&\PP\otb(\DD\otb\MM)\ar[r]^{\PP\otb\zeta^M} &\PP\otb(\MM\ota\CC)
}$$ what expands into the diagram
$$\xymatrix{
EGPBMA\ar[r]^{E\tau BMA}\ar[d]_\pr
&EPDBMA\ar[r]^{\sigma\chi MA}\ar[d]^{\pr}
&PBBDMA\ar[r]^{PBB\gamma MA}\ar[d]^\pr
&PBBMCA\ar[r]^{PB\alpha\psi}\ar[d]^\pr
&PBMAAC\ar[d]^\pr
\\
EGPB.MA\ar[dd]_\cong\ar[r]^{E\tau B.MA}
&EPDB.MA\ar[d]^{\sigma DB.MA}\ar@/_4pc/[dd]_{\cong}
&PB.BDMA\ar[d]^\cong\ar[r]^{PB.B\gamma A}
&PB.BMCA\ar[dd]^\cong\ar[r]^{PB\alpha CA}&PB.MAAC\ar[d]^\cong
\\
&PBDB.MA\ar[r]^\cong&PBDMA\ar[dr]^{PB\gamma A}&&PMAAC
\\
EGPMA\ar[r]^{E\tau MA}&EPDMA\ar[ru]^{\sigma DMA}\ar[r]^{EP\gamma
A} &EPMCA\ar[r]^{\sigma MCA}&PBMCA\ar[r]^{P\alpha
CA}&PMACA\ar[u]_{PMA\gamma} }$$

The only detail requiring explanation is the commutativity of the
hexagon below the second arrow in the upper row. To show that it
commutes one needs to expand it by inserting $PBDBMA$ in the
middle of the hexagon with a projection to $PBDB.MA$ and map
$PB\chi MA$ to $PBBDMA$ and also map $\sigma DBMA$ from the vertex
$EPDBMA$. Then the lower right corner of the split hexagon
commutes essentially by the compatibility of the unit $\eta^B$
with $\chi$.

Thus we obtained

\nxpoint {\bf Theorem.} {\it
The correspondences defined above, define a homomorphism
of bicategories (with the standard functoriality in pseudo-sense)
$$
\comc : \catbimod\longrightarrow\coring
$$
}

\section{Closing comments}

\nxpoint The operation of producing the lifting monad may be
also promoted to a canonical morphism of bicategories
from the bicategory of entwinings to the 2-category of monads
which act in the categories of right comodules
over variable coalgebras over variable rings. This is
in a complete analogy to one of the results in my earlier
article~\cite{skoda:ent}, so I will not bother writing
details here.

\nxpoint {\bf Acknowledgments.} I thank {\sc T. Brzezi\'nski} for
the encouragement and pointing out an important error in an early
version. The article has been written at IRB, Zagreb (partial
support of Croatia/MSES national projects 037-0372794-2807,
098-0000000-2865) and Max Planck Institute for Mathematics in Bonn
whom I thank for excellent working condition. My trips to Bonn
have been partly funded by DAAD/MSES bilateral project.

\footnotesize{

}


\begin{thebibliography}{99}
\bibitem{Semtriples}
{\sc H.~Appelgate, M.~Barr, J.~Beck, F.~W. Lawvere,
F.~E. J. Linton, E,~Manes, M.~Tierney, F.~Ulmer}, {\em Seminar on triples
and categorical homology theory}, ETH 1966/67, edited by B.~Eckmann, LNM 80,
Springer 1969.
\bibitem{Beck:distr}
{\sc Jon Beck}, {\em Distributive laws}, in~\cite{Semtriples},
119--140.
\bibitem{Bohm:internal}
{\sc G. B\"ohm}, {\em Internal bialgebroids, entwining structures
and corings}, AMS Contemp. Math. 376 (2005) 207-226; {\tt
arXiv:math.QA/0311244}
\bibitem{brz:bicat}
{\sc T. Brzezi\'nski, L. El Kaoutit, H. Gomez-Torecillas}, {\it
The bicategories of corings}, J.~Pure Appl Algebra 205: 510-541,
2006; {\tt math.RA/0408.5042}
\bibitem{BrzMajid:coalgBundles}
{\sc T. Brzezi\'nski, S. Majid}, {\it Coalgebra bundles}, Commun.
Math. Phys. 191:467-492, 1998; {\tt arXiv:q-alg/9602022}
\bibitem{brzWisbBk} {\sc T.~Brzezi\'nski, R. Wisbauer},
{\it Corings and comodules}, Cambridge Univ. Press 2003.
\bibitem{skoda:biact}
{\sc Z. \v{S}koda}, {\em Biactegories}, in preparation.
\bibitem{skoda:distr}
{\sc Z. \v{S}koda}, {\em Distributive laws for actions of monoidal
categories}, {\tt arXiv:math.CT/0406310}
\bibitem{skoda:ent}
{\sc Z. \v{S}koda}, {\em Equivariant monads and equivariant lifts
versus a 2-category of distributive laws}, {\tt arXiv:0707.1609}
\bibitem{Street:formal}
{\sc R. Street}, {\it The formal theory of monads}, JPAA 2, 149-168 (1972)
\end{thebibliography}
\end{document}